\documentclass[12pt]{amsart}

\input{amssymb.sty}

\newtheorem{letterthm}{Theorem}

\newtheorem{lettercly}[letterthm]{Corollary}
\newtheorem{theorem}{Theorem}[section]
\newtheorem{thm}[theorem]{Theorem}
\newtheorem{lem}[theorem]{Lemma}
\newtheorem{lemma}[theorem]{Lemma}
\newtheorem{cly}[theorem]{Corollary}
\newtheorem{prop}[theorem]{Proposition}
\newtheorem{proposition}[theorem]{Proposition}

\newtheorem{qu}[theorem]{Question}

\theoremstyle{definition}
\newtheorem{remark}[theorem]{Remark}

\def\Bbb{\mathbb}

\def\g{\gamma}
\def\G{\Gamma}

\def\Z{\Bbb{Z}}

\def\vc#1{\!\vec{\,#1}\,}
\def\S{\Sigma}

\def\isom{\cong}
\def\onto{{\kern3pt\to\kern-8pt\to\kern3pt}}

\def\<{\langle}
\def\>{\rangle}
\def\|{{\ |\ }}
\def\iff{{\leftrightarrow}}

\def\nm{\mbox{nm}}

\def\g{\gamma}
\def\G{\Gamma}

\catcode`\@=11
\def\serieslogo@{\relax}
\def\@setcopyright{\relax}
\catcode`\@=12

\begin{document}
\title[Subdirect products of groups]{Structure and finiteness properties\\
 of subdirect products of groups}

\author[Bridson]{Martin R.~Bridson}
\address{Martin~R.~Bridson \\
Mathematics Department \\
Huxley Building \\
Imperial College London \\
London SW7 2AZ \\
U.K. }
\email{m.bridson@imperial.ac.uk}

\author[Miller]{Charles~F. Miller~III }
\address{Charles~F. Miller~III\\
Department of Mathematics and Statistics\\
University of Melbourne\\
Melbourne 3010, Australia }
\email{ c.miller@ms.unimelb.edu.au }

\begin{abstract} We investigate the structure of subdirect products
of groups, particularly their finiteness properties. We pay special attention
to the subdirect products of free groups, surface groups 
and HNN extensions. We prove that a finitely presented subdirect product
of free and surface groups virtually contains a term of the lower central
series of the direct product or else fails to intersect one of the
direct summands.
This leads to a characterization of the finitely
presented subgroups of the direct product of 3 free or surface groups,
and to a solution to the conjugacy problem for arbitrary finitely
presented subgroups of direct products of surface groups.  
We obtain a
formula for the first homology of a subdirect product of two free groups
and use it to show there is no algorithm to determine the first homology
of a finitely generated subgroup. 
\end{abstract}

\subjclass{20F05, 20E07, 20F67, 20J05, 20F10}

\date{11 April 2007}

\keywords{subgroup structure, fibre product, subdirect product, free groups, surface groups,
homology of groups, HNN extension}

\maketitle

A useful structure theory for subgroups of finite direct products of 
groups has yet to be developed.  To begin to study such subgroups
it is natural to assume one knows about the subgroups of the direct factors and to concentrate on  subdirect products.
Recall that $G$ is termed a  {\em subdirect product} of the groups $A_1,\ldots, A_n$ if
$G \subseteq A_1\times\cdots \times A_n$ is a subgroup that projects 
surjectively to each factor.

Work by various authors  has exposed the surprisingly rich
structure to be found amongst the subdirect products of 
superficially-tame groups.   For example, in contrast to
the fact that
subdirect products  of abelian or nilpotent groups are
again in the specified class,
non-abelian free groups harbour a great diversity of 
subdirect products, including some with unsolvable 
decision problems \cite{cfm-thesis}. This diversity has long been
known, but it is only as a result of more recent work
by Baumslag-Roseblade \cite{gbjr} and Bridson-Howie-Miller-Short
\cite{bhms} that it has been understood as a phenomenon that is
intimately tied to the failure of various homological finiteness
conditions. 

From this background we pick-out the three strains of thought
to be developed in this article: the
subtlety of subdirect products in general; the usefulness of
finiteness properties in exploring this subtlety; and a special interest
in the subdirect products of free groups and their associates such as
surface groups. Further impetus for the study of subdirect products of
surface groups comes from the work of Delzant and Gromov \cite{DG}:
they proved that if a torsion-free group $\G$ is the fundamental
group of a compact K\"ahler manifold and $\G$ has sufficent multi-ended
splittings of the appropriate form,
 then there is a short exact sequence $1\to\Z^n\to\G_0\to S\to 1$,
where $S$ is a subdirect product of surface groups and $\G_0\subset\G$
is a subgroup of finite index.

\smallskip

Our first purpose in this article is to provide a systematic
and clarifying treatment of
the core issues that have emerged in the study of subgroups of direct products.
We focus on subdirect products of groups as objects worthy of study in their
own right. Following a discussion of some immediate consequences of the definiton
(Section \ref{s:general}), we address the  question of when such groups
are finitely generated and (more subtley) when they are finitely presented
(Section \ref{s:htpy}). We illustrate the general theory with a string of
explicit examples. In Section \ref{s:homol}
 we develop homological analogues of the results in Secton \ref{s:htpy}.
The main result in Section \ref{s:homol} provides one with a tool for analyzing
the second homology of subdirect products, a special case of which is the
following:

\begin{letterthm} Suppose that $G\leq F_1\times F_2$ is a subdirect product
of two free groups $F_1$ and $F_2$.   Let $L_i = G\cap F_i$. Then 
$$H_1(G,\Z) \cong  H_1(F_2,\Z) \oplus H_2(F_1/L_1,\Z)  \oplus C$$ 
where $C$ is a subgroup of $H_1(F_1,\Z)$ and hence is free abelian of
rank at most the rank of $F_1$.\qed
\end{letterthm}

In the light  of the Baumslag-Dyer-Miller construction \cite{gbedcm}, this yields:

\begin{lettercly} Let $F_1$ and $F_2$ be non-abelian free groups.  Then there are
continuously many subdirect products $G\leq F_1\times F_2$ 
having non-isomorphic $H_1(G,\Z)$.\qed
\end{lettercly}

Gordon \cite{cg} showed that $H_2$ of a finitely presented group is
not computable. By combining his construction with the above
theorem we deduce:

\begin{lettercly} If $F_1$ and $F_2$ are non-abelian free groups, 
there is no algorithm to compute for an arbitrary finitely generated subgroup 
$G\leq F_1\times F_2$ 
the (torsion-free) rank of $H_1(G,\Z)$, nor is there an algorithm 
to determine whether
$H_1(G,\Z)$ has non-zero torsion elements.  \qed
\end{lettercly}

In Section \ref{s:free} we sharpen our focus on subdirect products of 
free groups and surface groups and prove:

\begin{letterthm} \label{VirtNilp}
Let $\S_1,\ldots,\S_n$ be free groups  or surface groups and let 
$G\leq \S_1\times \cdots\times \S_n$ be a subdirect product
which intersects each factor $\S_i$ non-trivially.  
If $G$ is finitely presented, then each $\S_i$ contains a normal
subgroup $K_i$ of finite index such that 
$$\gamma_{n-1}(K_i) \subseteq G\cap \S_i \subseteq K_i.$$
Thus the quotients $\S_i/(G\cap \S_i)$ are virtually nilpotent 
of class at most $n-2$, and hence both $\S_i/(G\cap \S_i)$
and $G/(G\cap \S_i)$ are finitely presented.
\end{letterthm}

In the case of three factors, a 
refinement of the analysis used to prove Theorem \ref{VirtNilp} yields
the following characterization of finitely presented subgroups.

\begin{letterthm}\label{3fact} Let $\S_1,\S_2,\S_3$ be
finitely generated  free groups or surface groups 
and let  $G\leq \S_1\times \S_2 \times \S_3$ be a subdirect product
which intersects each factor $\S_i$ non-trivially.  
Then $G$ is finitely presented if and only if each $\S_i$ contains a normal
subgroup $K_i$ of finite index  such that the subgroup
 $G_0 = G \cap (K_1\times K_2\times K_3)$ satisfies the following condition:
 there is an abelian group 
$Q$ and epimorphisms $\varphi_i : K_i\to Q$ such that
$G_0\cap \S_i = \ker \varphi_i \ (i=1,2,3)$ and $G_0$ is the kernel of
the map $\varphi_1  + \varphi_2 +\varphi_3$.
\end{letterthm}

We are unsure what to expect for more than
three factors, but we have been unable to find any example
of a finitely presented subdirect product of free groups that 
intersects each of the direct factors but is
neither virtually a direct product of free groups nor virtually
normal with virtually abelian quotient, i.e.~{\em we do not know
if Theorem \ref{3fact} extends to an arbitrary (finite) number of factors}.
On the other hand,
even though Theorem \ref{VirtNilp} may not be optimal,
it is sufficient to provide considerable
control over the finitely presented subdirect products of free and surface
groups. For example, in Section \ref{decide} we use it to prove:

\begin{letterthm}\label{t:conj} If $D$ is a direct product of free and
surface groups, then every finitely presented subgroup $G\subset D$
has a solvable conjugacy problem. The membership problem
for $G$ is also solvable.
\end{letterthm}

In Section \ref{s:hnn} we shift our attention to subdirect products of
HNN extensions and establish a criterion for proving that certain
fibre products are not finitely presented. To exemplify
the utility of this criterion, we combine it
with an explicit calculation in group homology to prove the following
(cf.~\cite{bbhm}):

\begin{letterthm}\label{nfph}
 Let $ A = {\rm{BS}}(2,3) = \<b, t \| t^{-1}b^2t = b^3\>$
and let $q: A\to \mathbb Z = \<t\|\ \>$ be the map
defined by sending $t$ to $t$ and $b$ to 1.
Then the untwisted fibre product $G\subset A\times A$ 
associated to $q$ is finitely generated but not finitely presented,  and
$H_2(G,\Z)=0$.
\end{letterthm}

\section{\bf Generalities about subdirect products}\label{s:general}

If $G\leq A_1\times A_2$ is a subdirect product of groups
and we put $L_i = G\cap A_i$, then $G$ projects onto $A_2$ with kernel $L_1$.
The composition of this map with  $A_2\to A_2/L_2$
maps  $G$ onto
$A_2/L_2$ with kernel $L_1\times L_2$.  By
symmetry, 
 we have isomorphisms
$$ A_1/L_1 \cong G/(L_1\times L_2) \cong A_2/L_2.$$

Subdirect products of two groups are closely associated
to the fibre product construction in the category of groups. Recall
that,
associated to each pair of short exact sequences of
groups $1\to L_i \to A_i\overset{p_i}\to Q\to 1, \ i=1,2$, one has the 
{\em fibre product} or {\em pullback}
$$P=\{ (x,y) \in A_1\times A_2 \mid p_1(x)=p_2(y)\}.$$
Observe that  $L_1\times L_2 \leq P$ and that $P$ is generated
by  $L_1\times L_2$ together with any set of lifts $(t_1,t_2)$ of
a set of generators $p_i(t_i)$ for $Q$.  

It is clear that a fibre product is always a subdirect product.
Conversely, given a subdirect product $G\leq A_1\times A_2$,
 we can define $L_i := G\cap A_i$ and $Q := G/(L_1\times L_2)$, 
take $p_i$ to be the composition of the homomorphisms
$A_i\to A_i/L_i \overset{\sim}\to Q$, and regard $Q$
as ``the diagonal subgroup" in $A_1/L_1 \times A_2/L_2$. 
In more detail, given  $(x,y) \in A_1\times A_2$, we have
$(x,y)\in G$ if and only if $(x,y)(L_1\times L_2) \in Q$; that is, 
 if and only if $p_1(x) = (x,y)(L_1\times L_2) = p_2(y)$.
Thus $G$ is the fibre product of $p_1$ and $p_2$.

We summarize this discussion as follows:

\begin{prop}  A subgroup $G\leq A_1\times A_2$ 
is a subdirect product of $A_1$ and $A_2$ if and only if there
is a group $Q$ and surjections $p_i:A_i\to Q$ such that  $G$
is the fibre product of $p_1$ and $p_2$. \qed
\end{prop}

In the special case of a fibre product in which
$A_1 = A_2$, $L_1 = L_2$, and $p_1=p_2$, 
we shall call the fibre product {\em untwisted} (following \cite{bbhm}).

\begin{prop}\label{norm=vab}
Let $G\leq A_1\times\cdots\times A_n = D$ be a subdirect product of the groups 
$A_1,\ldots,A_n$ and let $L_i = G \cap A_i$.
Then the following are equivalent:
\begin{enumerate}
\item $G$ is normal in $D$;
\item each $A_i/L_i$ is abelian;
\item $G$ is the kernel of a homomorphism $\phi: D\to B$
where $B$ is an abelian group.
\end{enumerate}
\end{prop}

\proof First suppose that $G$ is normal in $D$.  Let $x,y\in A_i$.  Since
$G$ is a subdirect product there is an element of the form
$$\alpha  = (a_1,\ldots,a_{i-1},x,a_{i+1},\ldots,a_n)\in G.$$ 
Since $G$ is normal, it follows that 
$$\beta = (a_1,\ldots,a_{i-1},y^{-1}xy,a_{i+1},\ldots,a_n)\in G.$$
Hence $\alpha^{-1}\beta = (1,\ldots,1,x^{-1}y^{-1}xy,1\ldots,1) \in G$.
Thus $x^{-1}y^{-1}xy \in L_i = G\cap A_i$. Since $x$ and $y$ were
arbitrary elements of $A_i$, it follows that $A_i/L_i$ is abelian.
Now suppose each $A_i/L_i$ is abelian.  Then $G$ contains the derived group 
$[D,D]$ of $D$ and hence $G$ is the kernel of a homomorphism onto an abelian group.
That (3) implies (1) is trivial.  
\qed

\begin{cly}
Let $G\leq A_1\times A_2$ be a subdirect product of the groups $A_1$ and $A_2$
and let $L_i = G \cap A_i$.
Then $G$ is normal if and only if $G/(L_1\times L_2)$ 
(\ $\isom A_1/L_1 \isom A_2/L_2$) is abelian. \qed
\end{cly}

Observe that in this Corollary the group $G$ is the fibre product of two 
epimorphisms $p_1: A_1\to Q$ and $p_2: A_2\to Q$
where $Q=G/(L_1\times L_2)$.
Define $\phi: A_1\times A_2\to Q$ by $\phi(a_1,a_2) = p_1(a_1) p_2(a_2)^{-1}$.
Since $Q$ is abelian,  $\phi$ is a homomorphism.  Of course
$\phi(a_1,a_2) = 1$ if and only if $p_1(a_1) = p_2(a_2)$, so $G$ is the kernel of
$\phi$.  Clearly $\phi$ is surjective.  Thus 
$Q=G/(L_1\times L_2) \cong (A_1\times A_2)/G$. 

A normal subgroup $L$ of a group $A$ is said to be {\em finitely
normally generated} if it is the normal closure in $A$ of finitely many of 
its elements.  Observe that if $A$ is a finitely presented group,
then $L$ is finitely normally generated if and only if $A/L$
is finitely presented.

\section{\bf Generators and relations for a subdirect product}\label{s:htpy}
Suppose that $G \leq A\times B$ is a subdirect product of $A$ and $B$. 
We are interested in obtaining a presentation for $G$.
We define $L_A = G\cap A$ and $L_B=G\cap B$. Note that  $L_A$
and $L_B$ are normal subgroups of $G$.

Let $x_1=(a_1,b_1),x_2=(a_2,b_2),\ldots$ be a set of generators for $G$.
Then the $a_1,a_2,\ldots$ generate $A$ and the $b_1,b_2,\ldots$ generated $B$.
Choose presentations $A=\<a_1,a_2,\ldots \| r_1(a)=1, r_2(a)= 1,\ldots\>$
and $B=\<b_1,b_2,\ldots \| s_1(b)=1, s_2(b)= 1,\ldots\>$ on these generators.

If $s_i(x)$ is the word on the $x_i$ corresponding to $s_i(b)$, then
$s_i(x) = (s_i(a),s_i(b)) = (s_i(a),1)$ in $A\times B$ and hence
$s_i(x) \in L_A$.  More generally, if $w(x)$ is any word in the $x_i$
then $w(x)\in L_A$ if and only if $w(b) =_B 1$.  Hence $L_A$
is normally generated by the words $s_1(x),s_2(x),\ldots$.  Similarly
$L_B$ is normally generated by the words $r_1(x), r_2(x),\ldots$. 
Also $w(x)=_G 1$ if and only if $w(x)\in K_G = \nm(s_i(x))\cap \nm(r_j(x))$ 
where for instance $\nm(s_i(x))$ denotes the normal closure of the $s_i(x)$ 
in the free group on the $x_i$.  Thus $G$ can be presented
with generators $x_i$ and defining relations any set of normal generators for
$K_G$.  Notice that $K_G$ contains the commutator group 
$[\nm(s_i(x)),\nm(r_j(x))]$ and consequently 
$L_A$ and $L_B$ commute in $G$.

Of course if $G$ is finitely generated, then
$A$ and $B$ must be finitely generated. 
But the converse need not be true.

\subsection{Presenting finitely generated subdirect products}

As a consequence of the above discussion, we observe the following:

\begin{prop}  \label{fng} Suppose that $G$ is a subdirect product of two finitely
generated groups $A$ and $B$.  
\begin{enumerate}
\item If $G$ is finitely generated and $B$
is finitely presented, then  $G\cap A$ is finitely normally generated. 
\item If $G$ is finitely presented, then $B$
is finitely presented if and only if  $G\cap A$ is finitely 
normally generated. \qed
\end{enumerate} 
\end{prop}

Now suppose that both $G$ and $B$ are finitely presented on the generators
given above and that $G\cap A\neq 1$.  Unfortunately, it need not be true that
$A$ is finitely presented (see Example 5 below).

\subsection{Examples illustrating a diversity of finiteness properties}

We now give a series of examples which illustrate 
a variety of possibilities for finite generation and finite presentation of 
subdirect products which intersect both factors.

{\bf Example 1:}  It is a consequence of a theorem of Baumslag and Roseblade \cite{gbjr}
(see Theorem \ref{baumrose} below for a more complete statement), that a subdirect
product of two non-abelian free groups which intersects both factors is finitely 
presented if and only if it has finite index.  So ``most''  subdirect products of free
groups are not finitely presented. Trying to better understand this and related 
phenomena has been one of our motivations for this study.

As an explicit example, let $Q = \< t\|\>$ be infinite cyclic. Let $C=\< c_1,c_2\|\>$
be a free group mapping onto $Q$ by $c_1\mapsto t$ and $c_2\mapsto 1$.
Similarly let $D=\< d_1,d_2\|\>$ be free mapping onto $Q$
by $d_1\mapsto t$ and $d_2\mapsto 1$.  Let $G\leq C\times D$ be the
untwisted pullback of these two maps.  Then $G$ is generated by
$g_1=(c_1,d_1)$, $g_2=(c_2,d_2)$, $ g_3 = (c_2,1)$ and $g_4 = (1,d_2)$.
But $G$  has infinite index in $C\times D$ and so is not finitely presentable.
Notice that in this example $G\cap C$ and $G\cap D$ are both 
finitely normally generated.

\medskip

{\bf Example 2:} The second of our examples exhibits the nicest behaviour.
 Let $Q=\<c_1,c_2\| [c_1,c_2] = 1\>$, a free abelian group
of rank two.  Let $A=\<a_1,a_2\|\>$ be a free group 
and let
$$B=\<b_1,b_2\| [[b_1,b_2],b_1]=1,[[b_1,b_2],b_2]=1 \>$$ 
be free nilpotent group of class two and rank two.
Map $A$ and $B$ onto $Q$ via the maps $a_i\mapsto c_i$ and $b_i\mapsto c_i$.  
Let $G$ be the pullback of these two maps.  Then $G$ is generated
by the elements $x_1 = (a_1,b_1)$, $x_2 = (a_2,b_2)$, $x_3 =([a_1,a_2],1)$, 
$x_4 = (1,[b_1,b_2])$. 

To use the notation of the previous discussion, we re-present $A$ and $B$ on
the corresponding generators as 
\begin{eqnarray*}
A &=& \<a_1,a_2,a_3,a_4 \| a_4 = 1, [a_1,a_2]a_3^{-1}=1\>\\
B &=& \<b_1,b_2,b_3,b_4 \| [[b_1,b_2],b_1]=1,\\
&& \qquad  [[b_1,b_2],b_2]=1, b_3=1, [b_1,b_2]b_4^{-1}=1 \>.
\end{eqnarray*}
Then $L_A$ is normally
generated by the $s_i(x)$ of the previous discussion
which are $[[x_1,x_2],x_1]$, $[[x_1,x_2],x_1]$, $x_3$,
$[x_1,x_2]x_4^{-1}$.  
The subgroup $L_B$ is normally generated by the two elements 
$x_4$, $[x_1,x_2]x_3^{-1}$ (these are the $r_j(x)$).

In this example $G$ is finitely presented.  To see this observe that 
$A$ is free of rank two and $L_B$ is infinite cyclic with generator $[b_1,b_2]$, 
and $G$ is the split extension of $L_B$ by $A$. (In fact, since $L_B$
is central in $A\times B$, so $G\cong\mathbb Z\times F_2$.)

\smallskip

We record some general observations related to the salient points
of the preceding example. We remind the reader that a group $Q$ is
said to be of type $\mathcal F_3$ if there is an Eilenberg-Maclane
space $K(Q,1)$ with a finite 3-skeleton.

\begin{prop}  Let $A$ be a finitely generated free group, let
$p_1:A\to Q$ be an epimorphism, let $B$ be a finitely presented group
that fits into a short-exact sequence $1\to N \to B\overset{p_2}\to Q\to 1$,
and let $G\subset A\times B$ be the fibre product of $p_1$ and $p_2$. Then,
\begin{enumerate}
\item $G\cong A\ltimes N$;
\item if $N$ is finitely presented then $G$ is finitely presented;
\item if $N$ is finitely generated and $Q$ is of type $\mathcal F_3$
then $G$ is finitely presented.
\end{enumerate}
\end{prop}

The decomposition in (1) is obtained by splitting the projection of
$G$ to the first factor in $A\times B$, and in the light of (1)
assertion (2) is trivial. We omit the proof of (3), which is covered
by the arguments used in  \cite{bbms} to prove the 1-2-3 Theorem,
which provides a more subtle criterion for finite presentability of
fibre products.

\medskip

{\bf Example 3:} The essential feature of this example is that although
$A$ and $B$ are finitely presented, $G$ fails to be even finitely generated
as a result of $Q$ not being finitely presented.

Let $Q=\<c_1,c_2\| q_1(c),q_2(c),\ldots\>$ be a two
generator group which is not finitely presentable.  Let $A=\<a_1,a_2\|\>$
and $B=\<b_1,b_2\|\>$ be two free groups mapping onto $Q$ via the
maps $a_i\mapsto c_i$ and $b_i\mapsto c_i$.  Let $G$ be the
untwisted pullback of these two maps.  Then $G$ is generated
by the diagonal generators $x_1 = (a_1,b_1), x_2 = (a_2,b_2)$ together with
the elements $(q_1(a),1),(q_2(a),1),\ldots $. 
(Notice that $(1,q_i(b)) = q_i(x) (q_i(b),1)^{-1}$ so we do not need to add these
as generators.) 
Now since $Q$ is not finitely presented, no finite subset of these
generators suffice to generate $G$.  For if 
$(a_1,b_1)$, $(a_2,b_2))$, $(q_1(a),1))$, $\ldots$, $(q_n(a),1)$ generated $G$
then $Q \isom  A/L_A$ would be finitely
presented with presentation $\<a_1,a_2\| q_1(a)=1,\ldots, q_n(a)=1\>$.
Thus $G$ is not finitely generated.

\medskip

{\bf Example 4:} Let $Q$ be not finitely presentable as in Example 3  and let 
$$A = \< a_1,a_2,a_3 \| q_1(a) = 1, q_2(a) = 1,\ldots \>$$
where the words $q_i$ are as in the previous example. (We are using the
standard functional notation, so $a$ is the ordered alphabet $a_1,a_2,\dots$;
note that none of the relations $q_i(a)=1$  involve the
generator $a_3$.)
Similarly, let
$B = \< b_1,b_2,b_3 \| q_1(b) = 1, q_2(b) = 1,\ldots \>$, so
$B\isom A \isom Q * \<a_3\|\,\>$.  Map $A$ onto $Q$ by $a_1\mapsto c_1$,
$a_2\mapsto c_2$ and $a_3 \mapsto 1$, and  map $B$ onto $Q$ similarly. 
Since $Q$ is a free factor, neither $A$ nor $B$ is finitely presentable.

Again let  $G$ be the untwisted pullback of these two maps.  
In this example $G$ is finitely generated by the elements
$x_1 = (a_1,b_1)$, $ x_2 = (a_2,b_2),
x_3=(a_3,b_3)$ together with
 $ (a_3,1)$ and  $ (1,b_3)$ (one of these last two is redundant), but
 once again $G$ is not finitely presented. This time, the
 lack of finite presentability can be seen as a special case of
 Proposition \ref{fng}(2), for we are assuming $B$ is not finitely presentable
 and $L_A=G\cap A$ is normally generated by the single element $(a_3,1)$
 
 \medskip

{\bf Example 5:} Once again we take $Q$ as in Example 3 and
$$A = \< a_1,a_2,a_3 \| q_1(a) = 1, q_2(a) = 1,\ldots \>$$ as in example 3, but this time we
 let $B = \< b_1,b_2 \|\, \>$ be free of rank two.
We map $A$ onto $Q$ by $a_1\mapsto c_1$,
$a_2\mapsto c_2$ and $a_3 \mapsto 1$, and map $B$ onto $Q$ by $b_i\mapsto c_i$.
Again let  $G$ be the pullback of these two maps.  
In this example $G$ is generated by the finite collection of elements
$x_1 = (a_1,b_1)$, $ x_2 = (a_2,b_2)$, $x_3=(a_3,1)$. 
Again $L_A=G\cap A$ is the normal closure in
$A$ of $ x_3=(a_3,1)$.  Also $q_i(x) = (q_i(a),q_i(b)) = (1,q_i(b)) \in G$
and so $L_B$ is the normal closure in $B$ of the $(1,q_i(b))$.

In this example, it is more difficult to determine whether $G$ is finitely
presentable, but  in fact it is not since $L_A$ is not finitely generated
(see Theorem \ref{qjmthm} below).

\medskip

\subsection{Criteria for finite generation}

In our previous discussion we {\em assumed} we had a set of generators for 
the subdirect product $G\leq A\times B$.  We would like to know when $G$ 
actually has a {\em finite} set of generators.  Of course $A$ and $B$
must be finitely generated, so we assume this is the case.  Since $G$
is a subdirect product we can find finitely many elements 
$x_1=(a_1,b_1),\ldots,x_n=(a_n,b_n)$
in $G$ such that the $a_i$ generate $A$ and the $b_i$ generate $B$.
Denote by $H$ the subgroup generated by  $x_1,\ldots,x_n$. 
Then $H$ is a subgroup of $G$ but it may not contain all of either 
$L_A$ or $L_B$.  Observe that $H$ itself is a subdirect product of
$A$ and $B$. Also note that if $(u(a),v(b))\in G$ then 
$u(x)^{-1}(u(a),v(b)) = (1,u(b)^{-1}v(b)) \in L_B$ and 
$v(x)^{-1}(u(a),v(b)) = (v(a)^{-1}u(a),1) \in L_A$.
Thus $G = HL_A =HL_B$.  Hence if either $L_A$
 or $L_B$ is finitely normally generated, then $G$ will be finitely
generated. For instance, if $L_A$ is finitely normally generated
by $(z_1(a),1),\ldots,(z_n(a),1)$, then $G$ is generated by the $x_i$
and the $(z_j(a),1)$.  

Notice that in Example 4 the subgroup 
$L_A$ is finitely normally generated but $L_B$ is not.

We record the preceding general observation in the following proposition.

\begin{prop} Suppose that $G\leq A \times B$ is a subdirect product
of two finitely generated groups $A$ and  $B$. 
If  either $G\cap A$ or $G\cap B$ is 
finitely normally generated, then $G$ is finitely generated.  \qed
\end{prop}

Combining this with Proposition  \ref{fng},  we conclude the following:

\begin{cly} \label{fgnas}
Suppose that $G\leq A_1 \times A_2$ is the subdirect product
of two finitely presented groups $A_1$ and $A_2$.  Let $L_i = G \cap A_i$.
Then $G$ is finitely generated
if and only if  one (and hence both) of $A_1/L_1$ and $A_2/L_2$ are 
finitely presented. \qed
\end{cly}

\section{\bf Homological properties of subdirect products}\label{s:homol}

In this section we consider homological versions of the
 results from the previous section.  Recall that if
$A$ is finitely presented, then the integral homology groups $H_1(A,\Z)$
and $H_2(A,\Z)$ are both finitely generated.  On the other hand, there
exist non-finitely generated groups that have $H_1(A,\Z)=0$ 
as well as finitely generated groups $G$, 
with $H_2(G,\Z)$ finitely generated, which are not finitely presentable
(see, e.g., \cite{gbmult}, \cite{bbhm} or Theorem \ref{nfph} above). We also remind the reader
that there exist finitely presented groups whose
higher homology groups $H_n(G,\Z)$ are not finitely generated; explicit
examples due to Stallings 
and Bieri are described in the next section.

Consider a subdirect product $G\leq A\times B$ and let $L= G\cap A$.
As above we think of $G$ as an extension of $L$ by $B$. Conjugation
in $G$ and $A$ induces actions of these groups on their normal subgroup
$L$ and hence on $H_1(L,\Z)$. Since $G$ is a subdirect product,
its image in $\rm{Aut}(L)$ is the same as that of $A$ and hence 
$$H_0(G/L,H_1(L,\Z)) \cong H_0(A/L,H_1(L,\Z)). $$
 We make use of this observation to
prove the following:

\begin{thm}  \label{fgh1}
Let $A$ and $B$ be groups with both $H_1(-,\Z)$ and $H_2(-,\Z)$
finitely generated.  Suppose that $G\leq A\times B$ is a subdirect product of
$A$ and $B$. Then $H_1(G, \Z)$ is finitely generated if and only if one 
(and hence both) of $H_2(A/(G\cap A),\Z)$ and $H_2(B/(G\cap B),\Z)$ is
finitely generated.
\end{thm}

\proof Let $L= G\cap A$.
The usual five term exact sequence for $A/L$ gives the exactness of 
$$\cdots H_2(A,\Z) \to H_2(A/L,\Z)\to H_0(A/L,H_1(L,\Z))
 \to H_1(A,\Z) \cdots.$$
By hypothesis $H_1(A,\Z)$ and $H_2(A,\Z)$ are both finitely generated,
so $H_2(A/L, \Z))$ is finitely generated if and only if $H_0(A/L,H_1(L,\Z))$
is finitely generated.

Similarly, the five term exact sequence for $G/L$ gives the exactness of 
$$\cdots  H_2(G/L,\Z)\to H_0(G/L,H_1(L,\Z))
 \to H_1(G,\Z) \to H_1(G/L,\Z)\to 0.$$
By hypothesis $B = G/L$ and so $H_1(G/L,\Z)$ and $H_2(G/L,\Z)$ 
are both finitely generated.  Thus $H_1(G,\Z)$ is finitely generated 
if and only if $H_0(G/L,H_1(L,\Z))$ is finitely generated.

Since $H_0(G/L,H_1(L,\Z))$ and $H_0(A/L,H_1(L,\Z))$ are isomorphic, 
it follows that
$H_2(A/L, \Z))$ is finitely generated if and only if $H_1(G,\Z)$ 
is finitely generated, as claimed.  
By the symmetric argument, the other assertion follows.
\qed

\smallskip

If both factors $A$ and $B$ are free, the exact sequences
in the above proof yield more precise information, since in this case the
$H_1(-,\Z)$ are free abelian and $H_2(-,\Z) = 0$.  
We record this as follows:

\begin{cly} Suppose that $G\leq F_1\times F_2$ is a subdirect product
of two free groups $F_1$ and $F_2$.   Let $L_i = G\cap F_i$. Then 
$$H_1(G,\Z) \cong  H_1(F_2,\Z) \oplus H_2(F_1/L_1,\Z)  \oplus C$$ 
where $C= \ker (H_1(F_1,\Z)\onto H_1(F_1/L_1,\Z))$ and hence is free abelian of
rank at most the rank of $F_1$.\qed
\end{cly}

Since it is known how to construct two-generator groups with prescribed
countable $H_2(-,\Z)$ (see \cite{gbedcm}), one can apply the pull back 
construction to conclude the following result from \cite{gbjr}:

\begin{cly} Let $F_1$ and $F_2$ be non-abelian free groups.  Then there are
continuously many subdirect products $G\leq F_1\times F_2$ 
having non-isomorphic $H_1(G,\Z)$.\qed
\end{cly}

A theorem of Gordon asserts there is no algorithm to decide of a finitely presented
group $\G$ whether or not $H_2(\G,\Z) = 0$.  
This is proved (see \cite{cg} or \cite{cfm:msrisurvey}) by constructing a recursive
collection of finite presentations of groups $\G_i$
for which no such algorithm exists; these presentations
have a common finite set of symbols as generators, thus the
$\G_i$ come equipped with a surjection $F\to\G_i$ from a
fixed finitely generated free group.
It is easy to check from the construction that each of the groups 
involved is perfect.  Moreover one can easily arrange that 
each $\G_i$ is either the trivial group  or else
 $H_2(\G_i,\Z)\isom \Z \oplus (\Z/2\Z)$ 
(or any other fixed finitely generated abelian group).  
We now apply the above corollary to the pullback $G_i$ of two
copies of the presentation map $F \to \G_i$.  Since each $\G_i$ is perfect, 
we have 
$$H_1(G_i,\Z) \cong  H_1(F,\Z) \oplus H_2(\G_i,\Z)  \oplus H_1(F,\Z).$$ 
This proves the following result.

\begin{cly}   Let $F_1$ and $F_2$ be non-abelian free groups.  There is
a recursive collection of finitely generated subgroups $G_i$ of 
$F_1\times F_2$ such that there is no algorithm to compute the rank
of $H_1(G_i,\Z)$  or to determine whether it has any non-trivial
torsion elements. \qed
\end{cly}

As another application we note the following example.  
Let $F$ be a finitely generated
free group and suppose that $F/L$ has finitely generated  $H_2(F/L,\Z)$ but
is not finitely presented (cf.~Theorem \ref{nfph}).
Let $G\leq F\times F$ be the pullback or fibre product corresponding to 
this presentation. Then $H_1(G,\Z)$ is finitely generated
by Theorem \ref{fgh1}, but $G$ is not finitely generated by Corollary \ref{fgnas}.  
We record this as the following:

\begin{cly}  There is a subdirect product $G\leq F\times F$ of two finitely
generated free groups such that $H_1(G,\Z)$ is finitely generated
but $G$ is not finitely generated. \qed
\end{cly}

\section{\bf Subdirect products of free and surface groups}\label{s:free}

In this section and the next we focus on subdirect products
of free and surface groups.

\subsection{Background}
The results at the end of the previous section indicate how wild the
finitely generated subgroups of the direct product of two free groups
can be. But the following result of Baumslag and Roseblade shows that 
the only finitely presented subgroups are ``the obvious ones".

\begin{thm}[Baumslag and Roseblade \cite{gbjr}] \label{baumrose}
Let $F_1\times F_2$ be the
direct product of two free groups $F_1$ and $F_2$.  Suppose
that  $G \leq F_1\times F_2$ is a subgroup and define $L_i = G\cap F_i$.
\begin{enumerate}
\item If either $L_i = 1$ then $G$ is free.
\item If both $L_i$ are non-trivial and one of them is finitely generated, 
then $L_1\times L_2$ has finite index in $G$.
\item Otherwise, G is not finitely presented.
\end{enumerate}
\end{thm}

This result contains Grunewald's earlier result \cite{fg}
that in a direct product of two isomorphic
free groups, the untwisted fibre product $P\leq F\times F$ of a finite
presentation of an infinite group is finitely generated but not finitely
presented.  We remind the reader that such fibre products
can have unsolvable membership
problem and unsolvable conjugacy problem 
(see for instance \cite{cfm:msrisurvey}). 

\medskip

The following construction
shows that finitely presented subgroups of the direct product of more
than two free groups can be considerably more complicated than in the
case of two factors.

\medskip

{\bf Examples of Stallings and Bieri:} 
 Let $F_1=\<a_1,b_1\|\>$, $\ldots$,
$F_n=\<a_n,b_n\|\>$ be free groups of rank 2 and let $Q = \<c\|\>$ be an 
infinite cyclic group.  Let $\phi_n$ be the map from the direct product
$F_1\times \cdots\times F_n$ to $Q$ defined by $a_i\mapsto c$ and $b_i\mapsto c$.
Define ${\rm{SB}}_n = \ker \phi_n$.  It is easy to check that ${\rm{SB}}_n$ is a subdirect
product of the $F_i$ and that (for $n>1$) ${\rm{SB}}_n$ is finitely generated
by the elements $a_ib_i^{-1}$, $a_ia_j^{-1}$ and $b_ib_j^{-1}$.
Moreover, it can be shown \cite{bieri} that ${\rm{SB}}_n$ is of type 
${\rm{FP}}_{n-1}$ (even better,  type $\mathcal F_{n-1}$). But it is not of type ${\rm{FP}}_n$,
indeed $H_n({\rm{SB}}_n,\Z)$ is not finitely generated. (See \cite{stall1}
and \cite{bieri} for details.)

In order to relate these observations to our consideration of subdirect
products of two groups, we observe
 that projection onto the first $n-1$ factors maps ${\rm{SB}}_n$
surjectively onto $F_1\times \cdots\times F_{n-1}$ with kernel $L_n$
which is the normal closure in $F_n$ of $a_nb_n^{-1}$.  Further
 ${\rm{SB}}_{n-1} = {\rm{SB}}_n \cap (F_1\times \cdots\times F_{n-1})$.
Note that for $n>2$ the group ${\rm{SB}}_n$ is finitely presented.

\medskip

In the light of the diverse behaviour we have seen among the 
finitely generated subgroups of the direct product of two free groups, one
might expect that the above examples are just the first in a menagerie
of increasingly  exotic finitely presented subgroups in the case of 
three or more factors. However, somewhat to our surprise, this does
not appear to be the case. 

The first sign of  tameness among the finitely presented subgroups of
direct products of arbitrarily many free groups comes from the following
theorem of Bridson, Howie, Miller and Short \cite{bhms}, which shows that
whatever wildness exists may be detected at the level of homology.

\begin{thm}[\cite{bhms}]
Let $F_1,\dots,F_n$ be free groups. A subgroup
 $G\le F_1\times\dots\times F_n$ is of type ${\rm{FP}}_n$ 
 if and only if it has a subgroup of finite index that
is itself a direct product of (at most $n$) finitely generated free groups.
\end{thm}

Thus, on the one hand, we know that the only homologically-tame subgroups
of a direct product of free groups are the obvious ones. On the other
hand, we have a specific method for constructing  examples of homologically-wild
subgroups coming from the construction of Stallings and Bieri. Moreover one
has essentially complete knowledge of the latter situation, because the BNS
invariants of direct products of free groups have been calculated
\cite{meinert}, providing a complete classification of the finiteness
properties of the kernels of maps from $F_1\times\dots\times F_n$ to abelian
groups.

Our repeated failure to construct finitely-presented subgroups that
that are neither ${\rm{FP}}_\infty$ nor of Stallings-Bieri type leads us to pose
the following:

\begin{qu}\label{quest} Let  $D = F_1\times \cdots\times F_n$
be a direct product of free groups (of various ranks) and 
let $p_i:G\to F_i$ be the natural projection. Let
$G\subset D$ be a subgroup that intersects each $F_i$ non-trivially.

If $G$ is finitely presented but not of type ${\rm{FP}}_n$, then does $G$ have
a subgroup of finite index $G_0$ which is normal in 
$p_1(G_0)\times \cdots \times p_n(G_0)$ with abelian quotient?
\end{qu}

Theorems \ref{virtnilp} and \ref{3fact} add
considerable interest to this question.  In the course
of proving these results we will make use of the
 following theorem  from \cite{cfm:dpwithfree},
which was used there to give a straightforward proof of the 
Baumslag-Roseblade Theorem:

\begin{thm}[Miller\cite{cfm:dpwithfree}]    \label{qjmthm}
Let  $A\times F$ be the direct product of a group $A$ with a free group $F$.
Suppose that $G \leq A\times F$ is a subgroup which intersects $F$ non-trivially.
If $G$ is finitely presented, then $L = G\cap A$ is finitely generated. \qed
\end{thm}

\subsection{Preparatory results for the surface case}
Since we want our results to apply to both free and surface groups,
we must also prove the analog of the above result for compact surfaces.  One ingredient
of the proof is the well known fact \cite{scott} that a 
finitely generated, non-trivial normal
subgroup of a compact surface group of negative Euler characteristic
must have finite index.   A second important ingredient is the following
substitute for the use made in 
\cite{cfm:dpwithfree} of a theorem of M. Hall.

\begin{lem} \label{mhall-surface}
Let $\Pi$ be the fundamental group of a closed surface $S$ and let $H\le \Pi$
be a non-cyclic 2-generator subgroup. Then there exist
elements $a_1,b_1\in H$ that, in a finite index subgroup $\Pi_0\le\Pi$, serve as the
beginning generators in a standard presentation 
$$\Pi_0=\<a_1,b_1,\ldots,a_g,b_g\mid [a_1,b_1]\cdots [a_g,b_g] = 1\>.$$
Equivalently, there a finite-sheeted covering $\hat S$ of $\S$ such that
$\<a_1,b_1\>\cap\pi_1\hat S$  is the fundamental group of a subsurface
of positive genus.
\end{lem}

\begin{proof}
By a theorem of Scott \cite{scott}, 
given a finitely generated subgroup $H$ of a surface group,
one may pass to a finite-sheeted cover $\hat S$ so that $H$ is the fundamental
group of a subsurface $T$ onto which $\hat S$ retracts. If
$H$ is a 2-generator free group, then this surface is either a 
once-punctured torus or a thrice-punctured sphere. In the former
case, we are done. In the latter case, there is a 4-sheeted
cover $\hat T$ of $T$ that is a 4-punctured torus;
let $c_1,c_2,c_3,c_4\in\pi_1T$ be coset representatives
of $\pi_1\hat T$. Using the
fact that surface groups are subgroup separable, we find a finite-sheeted
cover of $\hat S$ to which none of $c_1,c_2,c_3,c_4$ lift. A component
of the preimage of $T$ in this covering is a cover of $\hat T$ and
hence has positive genus.
\end{proof}

An alternative 
proof of this lemma can be derived from a theorem  of 
Bridson and Howie (Corollary 3.2  of \cite{mbjh2}).
The type of ``positive-genus" argument we used
has been exploited 
in \cite{mbmthw} to obtain results about elementarily free groups.
 
 We are now ready 
to prove the analog of Theorem \ref{qjmthm} for surface groups.

\begin{thm} \label{qjm-surface} 
Let $\S$ be the fundamental group of a compact surface
other than the Klein bottle and torus, and let $A$ be an arbitrary group.
Let $G \leq A\times \S$ be a subgroup that intersects $\S$ non-trivially.
If $G$ is finitely presented, then $G\cap A$ is finitely generated. 
\end{thm}

\begin{proof} Let $L=G\cap A$ and
let $p:G\to \S$ be the standard projection. If the surface
is a sphere or projective plane, the statement is trivial. 
More generally, if $G\cap\S$ has finite index in $p(G)$, then $G$ contains
$(G\cap A)\times (G\cap \S)$ as a subgroup of finite index.
And since $L$ is a retract of this subgroup, it is finitely
presented.

The finitely generated group $p(G)\subset\S$ is either free, in which
case we are done by theorem \ref{qjmthm}, or else it is again the
fundamental group of a closed surface. Thus there is no loss
of generality in assuming that   $p(G)=\S$.
This forces $G\cap\S$ to be normal in $\S$.
Thus we are reduced to the case where $G\cap \S\subset\S$
is a non-trivial normal subgroup of infinite index. 
It follows from Lemma \ref{mhall-surface} that by replacing
$\S$ with a subgroup of finite index  and taking the preimage in $G$,
we may assume that $\S$ has a presenetation of the form 
$$\<a_1,b_1,\ldots,a_g,b_g\mid [a_1,b_1]\cdots [a_g,b_g] = 1\>$$
where $a_1$ and $b_1$ both lie in $G\cap \S$.   The defining
relation is equivalent to the equation 
$$a_1^{-1} b_1 a_1 = b_1 [a_2,b_2]\cdots [a_g,b_g]$$
and so we regard $\S$ as the HNN extension with stable letter
$a_1$ which conjugates the cyclic subgroup $\<b_1\>$ to the 
cyclic subgroup generated by the right hand side.  
Note that $b_1,a_2,b_2,\ldots,b_g$ freely generate a free subgroup
of $\S$.   

We proceed as in \cite{cfm:dpwithfree}.
For each $i=2,\ldots,g$ pick a lift $\hat{a}_i\in p^{-1}(a_i)$ and
 $\hat{b}_i\in p^{-1}(b_i)$ in $G$ .  Observe that 
 $$   [a_2,b_2]\cdots [a_g,b_g]  =_G
[\hat{a}_2,\hat{b}_2]\cdots [\hat{a}_g,\hat{b}_g] \cdot c_1$$
 for some element $c_1\in L=G\cap A$.
 
Since $b_1,\hat{a}_2,\hat{b}_2,\ldots,\hat{b}_g$ are the pre-image
of a free basis, they  freely generate a free subgroup of $G$.
Hence the subgroup $H$ of $G$ generated by $L$ together with
these elements has the structure of an 
HNN extension of $L$.
The associated subgroup for each stable letter is $L$ itself and
we note that $b_1$ acts trivially on $L$ since $b_1\in G\cap\S$.

Now  $G$ is an extension of
 $H$ with stable letter $a_1$  and is finitely presented, 
 so it can be generated by
 $ a_1,b_1,\hat{a}_2,\ldots,\hat{b}_g$ together with finitely many
 elements $c_1,\ldots,c_n\in L$ (including the previously chosen
 $c_1$).  Thus we can present $G$ as
 $$G=\<  a_1,b_1,\hat{a}_2,\ldots,\hat{b}_g, c_1,\ldots,c_n \mid
  a_1^{-1}b_1a_1 = [\hat{a}_2,\hat{b}_2]\cdots [\hat{a}_g,\hat{b}_g] \cdot c_1,
$$
 $$\hfill\qquad 
  a_1^{-1} d a_1 = d \ (d\in L), 
 \mbox{ relations of } H\>.$$
 The associated subgroup for the stable letter $a_1$ is $L\times\<b_1\>$.
 Observe that $H$ is generated by the generators other than $a_1$ since
 the action of $a_1$ on $L$ is trivial.
 Hence, because $G$ is finitely presented,  $L\times\<b_1\>$ is finitely generated,
 and so $L$ is finitely generated. This completes the proof.
\end{proof}

\subsection{The Virtually-Nilpotent Quotients Theorem}

The following theorem controls the way in which
a finitely presented subdirect product of free 
and surface groups can intersects the direct factors
of the ambient group.

The reader will recall that
the $m$-th term of the lower central series of
a group $H$ is defined inductively by $\gamma_1(H)=H$ and 
$\gamma_{m+1}(H)=[\gamma_m(H),H]$. And
$H$ is defined to be nilpotent of class $c$ if $\gamma_{c+1}(H)=1$.

\begin{thm} \label{virtnilp}
Let $\S_1,\ldots,\S_n$ be free or surface groups and let 
$G\leq D = \S_1\times \cdots\times \S_n$ be a subdirect product
which intersects each factor $\S_i$ non-trivially.  
If $G$ is finitely presented, then each $\S_i$ contains a normal
subgroup $K_i$ of finite index such that 
$$\gamma_{n-1}(K_i) \subseteq G\cap \S_i \subseteq K_i.$$
Thus the quotients $\S_i/(G\cap \S_i)$ and
$$D/((G\cap\S_1) \times \cdots \times (G\cap\S_n))$$
are virtually nilpotent 
of class at most $n-2$. Hence both $\S_i/(G\cap \S_i)$
and $G/(G\cap \S_i)$ are finitely presented, and consequently
the projection of $G$ into the product
of any $j<n$ factors is finitely presented.
\end{thm}

Note that in the case $n=2$ the conclusion is that $K_i = G\cap \S_i$
which implies the Baumslag-Roseblade Theorem.  In case $n=3$,
the conclusion is that $[K_i,K_i]\subseteq G\cap \S_i \subseteq K_i$
so that the $\S_i/(G\cap \S_i)$ are virtually abelian, as happens
for the  Stallings-Bieri examples ${\rm{SB}}_n$. 

\begin{proof}
Let $p_i:G\to \S_i$ be the induced projection maps 
and put $N_i = \ker p_i$
and $L_i = G\cap \S_i$. Since $G$ is finitely presented, 
by Theorem \ref{qjmthm}
each $N_i$ is finitely generated.  Because $G$ is subdirect, 
each $\S_i$
is finitely generated, and the $L_i$ are normal in $\S_i$ as well as in $G$.

Again since $G$ is subdirect, for $j\neq i$, the projection $p_j(N_i)$ is 
normal in $\S_j$.  Now $L_j\subseteq p_j(N_i)$ so $p_j(N_i)$ is a non-trivial
finitely generated normal subgroup of the free group $\S_j$, and hence has
finite index in $\S_j$.

For notational simplicity we focus on the case $j=1$ and note that similar
arguments work for the remaining $j=2,\ldots,n$.  
Define $$K_1 = p_1(N_2)\cap \cdots \cap p_1(N_n).$$
We note that $L_1 = N_2\cap \cdots \cap N_n \subseteq K_1$.
For any choice of $n-1$ elements $x_2,\ldots,x_n \in K_1$, there
are elements $y_i \in N_i$ with $p_1(y_i) = x_i$ for $i=2,\ldots,n$.
Observe for example that $y_2$ has the form $(x_2,1,z_{2,3},\ldots,z_{2,n})$
and  $y_3$ has the form $(x_3,z_{3,2},1,\ldots,z_{3,n})$.  
Hence their commutator is
$$[y_2,y_3] = ([x_2,x_3],1,1,[z_{2,4},z_{3,4}],\ldots,[z_{2,n},z_{3,n}]) \in N_2\cap N_3.$$
On forming an $(n-1)$-fold commutator such as $[y_2,y_3,\ldots,y_n]$ one obtains
$$[y_2,y_3,\ldots,y_n] = ([x_2,x_3,\ldots,x_n],1,1,\ldots,1)\in G\cap \S_1 =  L_1,$$
and similarly for other commutator arrangements.  
Hence $\gamma_{n-1}(K_1) \subseteq L_1$.  But we know $L_1\subseteq K_1$, so
 $$\gamma_{n-1}(K_1) \subseteq L_1 \subseteq K_1$$ as desired.

The remaining assertions follow from the fact that 
finitely generated nilpotent groups are finitely presented, since this
implies that
any normal subgroup of a finitely generated free group that contains 
a term of the lower central series is finitely
normally generated.  This completes the proof.
\end{proof}

\medskip

In case $n=3$ we can actually characterize the finitely presented subdirect
products of $\S_1\times \S_2\times \S_3$.   

\begin{lem} \label{onto-pairs}
Let $\S_1,\S_2,\S_3$ be free groups or surface groups and let 
$G\leq \S_1\times \S_2 \times \S_3$ be a subdirect product
which intersects each factor $\S_i$ non-trivially.  
If $G$ is finitely presented, then each $\S_i$ contains a normal
subgroup $K_i$ of finite index such that the projections
of $G_0 = G \cap (K_1\times K_2\times K_3)$  to
pairs of factors,
$$p_{i j}: G_0 \to K_i \times  K_j$$
 $i<j$, are surjective.  
In particular the projections of $G_0$ to the $K_i$
are also surjective.
\end{lem}

\begin{proof}  Since the hypotheses are the same, we may 
continue with the notation of the proof of the previous theorem.
Let $k_1\in K_1$.  By the definition of $K_1$ there are triples
$(k_1,1,y_2)\in N_2\subset G$ and $(k_1,x_3,1)\in N_3\subset G$.  Thus
$(1,x_3^{-1},y_2)\in N_1$.  Therefore $x_3\in p_2(N_1)\cap p_2(N_3) = K_2$
and $y_2\in p_3(N_1)\cap p_3(N_2) = K_3$.  

It follows that $(k_1,1,y_2) \in N_2\cap (K_1\times K_2\times K_3)$ and
hence $(k_1,1)\in p_{12}(G_0)$. Also 
$(k_1,x_3,1) \in N_3\cap (K_1\times K_2\times K_3)$
and so  $(k_1,1)\in p_{13}(G_0)$.   Similar calculations apply
for other factors $K_i$ and so $G_0$ projects onto pairs as claimed.
\end{proof}

 Here is a lemma about 3 factors
with projections onto any 2 factors which works for any
groups.

\begin{lem} Let
$G\leq A_1\times A_2  \times A_3$ be a subdirect product
which projects surjectively onto any product of two factors
and let $L_i = G\cap A_i$.  Then there is an abelian group 
$Q$ and epimorphisms $\varphi_i : A_i\to Q$ such that
$L_i = \ker \varphi_i \ (i=1,2,3)$ and $G$ is the kernel of
the map $\varphi_1  + \varphi_2 +\varphi_3$.
\end{lem}

\begin{proof}  Since the projection to $A_2\times A_3$ is surjective,
we can think of $G$ as a subdirect product of the two groups 
$A_1 $ and $ ( A_2\times A_3)$.  Hence $G/L_1 \isom G/( G\cap (A_2\times A_3))$.
If we put $Q = G/L_1$ then $G$ can also be viewed as the pullback 
(fibre product) of the quotient map $\varphi_1: G\to Q$ and a surjection $\psi: A_2\times A_3 \to Q$.

Now $\psi$ restricts to homomorphisms $\theta_2:A_2\to Q$ and 
$\theta_3:A_3 \to Q$ with $\psi((a_2,a_3)) = \theta(a_2)\cdot \theta(a_3)$.
Since the projection from $G$ to  $A_1\times A_2$ is surjective,
for any $x\in A_1$ there is a $z\in A_3$ such that $(x,1,z)\in G$.
Since $G$ is the pullback of $\varphi_1$ and $\psi$, we have
$\varphi_1(x) = \psi((1,z))= \theta_3(z)$.    Because $x\in A_1$ was
arbitrary, it follows that $\theta_3$ has the same image as $\varphi_1$
and hence maps $A_3$ onto $Q$ with kernel $L_3$.

Similarly since  the projection from $G$ to  $A_1\times A_3$ is surjective,
$\theta_2$ maps $A_2$ surjectively onto $Q$ with kernel $L_2$.
But the images of $A_2$ and $A_3$ commute in $Q$ and so $Q$
must be abelian.  Changing to additive notation and defining
 $\varphi_2 = - \theta_2$ and $\varphi_3 = - \theta_3$, it follows
 that $G$ is the kernel of the map kernel of
the map  $\varphi_1  + \varphi_2 +\varphi_3$.  This completes the proof.
\end{proof}

Meinert \cite{meinert} has calculated the Bieri-Neumann-Strebel
invariants for direct products $D$ of finitely many
finitely generated free groups. In
particular he has calculated which homomorphisms from
$D$ to an abelian group have a finitely presented kernel.
By combining his result with the preceding two lemmas we obtain 
a characterization of finitely presented subdirect products of 3 free or
surface groups.

\begin{thm} \label{3fpchar}
Let $\S_1,\S_2,\S_3$ be finitely generated free groups or surface
groups and let 
$G\leq \S_1\times \S_2 \times \S_3$ be a subdirect product
which intersects each factor $\S_i$ non-trivially.  
Then $G$ is finitely presented if and only if  $\S_i$ contains a normal
subgroup $K_i$ of finite index  such that the subgroup
 $G_0 = G \cap (K_1\times K_2\times K_3)$ satisfying the following condition:
 there is an abelian group 
$Q$ and epimorphisms $\varphi_i : \S_i\to Q$ such that
$G_0\cap \S_i = \ker \varphi_i \ (i=1,2,3)$ and $G_0$ is the kernel of
the map $\varphi_1  + \varphi_2 +\varphi_3$.
\end{thm}

\begin{proof} The necessity of the given condition 
is established by the preceding two lemmas and its sufficiency
in the case of free groups is a special case of Meinert's
theorem \cite{meinert}. Thus we need only argue that
sufficiency in the case of surface groups follows from Meinert's
result.

To this end, for $i=1,2,3$ we choose an epimorphism
$\pi_i:F_i \to \S_i$ where $F_i$ is a finitely generated
 free group and the kernel is either trivial or the
 normal closure of a single product $c_i$ of
 commutators (i.e.~the standard surface relation).
 For notational convenience we define $c_i=1$ if
 $\S_i$ is free.
 
 Consider 
 the composition $\Phi_i:=\phi_i\pi_i:F_i\to Q$.
 Meinert's theorem tells us that
 the kernel $\Gamma_0\subset F_1\times F_2\times
 F_3$ of $\Phi_1+\Phi_2+\Phi_3$ is finitely presented.
 
 The kernel $I$ of the
 map $F_1\times F_2\times F_3\to\S_1\times
 \S_2\times\S_3$ induced by the $\pi_i$ is the normal
 closure of $C=\{c_1,c_2,c_3\}$, and $G_0=\G_0/I$.
 Since $\G_0$ is subdirect, the normal closure of $C$
 in $\G_0$ is the same as its normal closure in 
 $F_1\times F_2\times F_3$. It follows that a presentation
 of $G_0$ can be obtained by adding just three relations
 to a presentation for $\Gamma_0$. In particular, 
 $G_0$ is finitely presented.
 \end{proof}
 
Thus Question \ref{quest} has a positive answer for $n=3$.  We have on occasion
thought there is a positive answer for $n=4$, but a satisfactory proof has yet
to emerge.  For $n>4$ we are unsure of what to expect.

\begin{remark}
The restrictions that we have obtained concerning 
subdirect products of surface groups do not extend
to subdirect products of arbitrary hyperbolic groups, but it
appears that they do extend to subdirect products of
fully residually free groups. The same is true of the results
presented in the next section. We shall explore these matters
in a future article with J.~Howie and H.~Short.
\end{remark}

\section{\bf Decision problems for finitely presented subgroups}\label{decide}

The  results in this section should be contrasted with the fact
that if $\S_1$ and $\S_2$ are non-abelian free or surface
groups, then there are finitely generated subgroups $H\subset\S_1\times\S_2$
for which the conjugacy problem and membership problem are
unsolvable \cite{cfm-thesis}.

\subsection{The conjugacy problem}
We shall prove that finitely presented subgroups of
direct products of surface groups have a solvable conjugacy
problem. 
Our proof relies heavily on
the structure of such subgroups as described in 
Theorem \ref{virtnilp}. With this structure in
hand, we can adapt the argument used in the
proof of  Theorem 3.1
of \cite{mb-haef},
where it was proved that if $G$ is a bicombable group,
$N\subset G$ is a normal subgroup and the generalized word problem for $G/N$ is solvable, then $N$ has a solvable conjugacy problem. 

The class of bicombable groups contains the hyperbolic groups and
is closed under finite direct products and
the passage to subgroups of finite index. It follows that 
subgroups of finite index in direct products of free and surface groups are bicombable. The properties of bicombable
groups that we need are all classical and easy to prove in
the case of such subgroups, but we retain the greater
generality with an eye to future applications.

\smallskip
 With  some effort, one can give
an elementary proof of the
following lemma using induction on the
nilpotency class. A more elegant argument due to 
Lo (Algorithm 6.1 of \cite{Lo})  provides an algorithm that is practical
for computer implementation.

\begin{lemma}\label{Lo} If $Q$ is a finitely generated
nilpotent group, then there is an
algorithm that, given finite sets $S,T\subset Q$ and $q\in Q$,
will decide if $q\<S\>$ intersects $\<T\>$ non-trivially. \qed
\end{lemma}

The adaptation of Theorem 3.1 of \cite{mb-haef} that
we need in the present context is the following.

\begin{proposition}\label{solvC} 
Let $\G$ be a bicombable group, let $H\subset \G$ be a 
subgroup, and suppose that there exists
a subgroup $L\subset H$ normal in $\G$ such that 
$\G/L$ is nilpotent. Then $H$ has a solvable conjugacy
problem.
\end{proposition}

\begin{proof} The properties of a bicombable group $\G$ that
we need here are (1) the conjugacy problem is solvable in $\G$ and
(2) there is an algorithm that, given $g\in \G$
as a word in the generators of $G$, will
calculate a finite generating set for the centralizer of $g$.
(The second fact has its origins in the work of
Gersten and Short \cite{GS}; the running time of the algorithm
depends on the length of the word representing $g$ and on the
fellow-traveller constant of the bicombing.) The reader should have
little difficulty in supplying their own proof of these
facts in the case where
$\Gamma$ is a product of free and surface groups.

Given $x,y\in H$ (as words in the generators of $\G$)
we use the positive solution to the conjugacy problem in
$\G$ to determine if there exists $\gamma\in\G$ such that 
$\g x\g^{-1}=y$. If no
such $\g$ exists, we stop and declare that $x,y$ are not conjugate in $H$.
If $\g$ does exist then we find it and consider
$$\g C=\{g\in \Gamma\mid gxg^{-1}=y\},$$
where $C$ is the centralizer of $\g$ in $\G$.
Note that $x$ is conjugate to $y$ in $H$ if and only
if $\g C\cap H$ is non-empty.

We employ the algorithm from (2) to compute a finite generating
set $\hat S$ for $C$. We then employ Lo's algorithm
(Lemma \ref{Lo}) in the nilpotent group $\G/L$ 
to determine if the image of $\gamma C$ intersects the
image of $H$. Since $L\subset H$, this intersection is
non-trivial (and hence $x$ is conjugate to $y$) if and
only if $\g C\cap H$ is non-empty.
\end{proof}

We need a further lemma in order to make full use of the
preceding proposition. We remind
the reader that
the solvability of the conjugacy problem
does not in general pass to subgroups or overgroups of finite index \cite{CM}.

A group $G$ is said to have {\em unique roots} if for all
$x,y\in G$ and $n\neq 0$ one has $x=y\ \iff\  x^n=y^n$. Torsion-free
hyperbolic groups and their direct products have this property.

\begin{lemma}\label{findex} Suppose $G$ is a group in which roots are
unique and  $H\subset G$ is a subgroup of finite index.
If the conjugacy problem for $H$  is solvable, then
the conjugacy problem for $G$ is solvable.
\end{lemma}

\begin{proof} Let $m_0$ be the index of $H$ in $G$ and
let $m=m_0!$. Given $x,y,g\in G$, since roots are unique
$x^m=gy^mg^{-1}$ if and only if $x=gyg^{-1}$.
Thus $x,y$ are conjugate  in $G$ if and only if
$x^m,y^m$  are conjugate in $G$. Note $x^m,y^m\in H$.

If $c_1,\dots,c_{m_0}$
are coset representatives for $H$ in $G$ and $x_i:=c_ixc_i^{-1}$,
then $x^m$ is conjugate to $y^m$ in $G$ if and only if 
$y^m$ is conjugate to at least one of $x_i^m$ in $H$.

Combining these two observations, we see that deciding if
$x$ is conjugate to $y$ in $G$ reduces to deciding if
one of finitely many conjugacy relations holds in $H$. This
completes the proof.
\end{proof}

\begin{thm}  If $D$ is a direct product of free and
surface groups, then every finitely presented subgroup of $D$
 has a solvable conjugacy problem.
\end{thm}

\begin{proof}
Projecting $D$ away from direct factors that intersect $G$ trivially, and
replacing each of the remaining factors by the projection of $G$ to that
factor, we see that there is no loss of generality in assuming that $G$
is a subdirect product of
$D=\S_1\times\dots\times\S_n$ and 
that each $L_i=\S_i\cap G$ is non-trivial.

Theorem \ref{virtnilp} tells us that 
$L=L_1\times\dots\times L_n$
is normal in $D$ and $D/L$ is virtually nilpotent. Let
$N$ be a nilpotent subgroup of finite index in $D/L$,
let $D_0$ be its inverse image in $D$ and let $G_0=D_0\cap G$.

We are now in the situation of Proposition \ref{solvC} with
$\G=D_0$ and $H=G_0$. Thus $G_0$ has a solvable conjugacy problem.

Finally, since roots are unique in surface groups, they are
unique in $D$. Therefore Lemma \ref{findex} applies and
we conclude that the conjugacy problem for $G$ is solvable.
\end{proof}
 
\subsection{The membership problem}

In the course of proving our next theorem we will need the
following technical observation.

\begin{lemma}\label{l:makeP}
If $\Sigma$ is a finitely generated
free or surface group, then there
is an algorithm that, given a finite set $X\subset\Sigma$, will
output a finite presentation for the subgroup generated by $X$.
\end{lemma}

\begin{proof} Let $G$ be the subgroup generated by $X$.
The lemma is a simple consequence of the fact that $\Sigma$
has a subgroup of finite index that retracts onto $G$. This fact
is due to M.~Hall \cite{mh} in the case of free groups
and P.~Scott \cite{scott} in the case of surface
groups.

In more detail, running through the finite-index subgroups
$\Sigma_0\subset\Sigma$, one calculates a presentation for
$\Sigma_0$ and one attempts to express the elements $x\in X$
as words $u_x$ in the generators of $\Sigma_0$ by listing all
words in these generators and using the solution to the word problem
in $\Sigma$ to check equality. When words $u_x$ have been
found for all $x\in X$, one begins a naive search for homomorphisms
$\phi:\Sigma_0\to G$: products $v_b$
of the letters $X^{\pm 1}$ are chosen
as putative  images
for the generators $b$ of $\Sigma_0$ and the solution
to the word problem in $\Sigma$ is used to check if the defining
relations of $\Sigma_0$ are respected by this choice. When a 
homomorphism $\phi$ is found, one again uses the solution to the
word problem in $\Sigma$ to check if $\phi(u_x)=x$ for all $x\in X$.

One applies this procedure to all $\Sigma_0$, proceeding
in a diagonal manner. The theorems of Hall and Scott assure us
that it will eventually terminate, at which point we have
a presentation $\Sigma_0 = \langle B\mid R\rangle$ and 
words $\{u_x,\, v_b\mid x\in X,\, b\in B\}$. The desired presentation
of $G$ is then 
$\langle B\cup X \mid R,\, xu_x^{-1},\, bv_b^{-1}\ (x\in X,\,b\in B)\rangle$.
\end{proof}

\begin{thm}\label{t:memb} If $D$  
is the direct product
of finitely many finitely generated
free and surface groups and $G\subset D$ is a 
finitely presented subgroup, then the membership problem for 
$G$ is decidable, i.e.~there is an algorithm which, 
given $h\in D$ (as a word in the generators)  will
determine whether or not $h\in G$.
\end{thm}

\begin{proof}  We proceed by induction on the
number of factors in $D=\S_1\times\dots\times\S_n$.
For $n=1$ the assertion of the theorem is well-known,
in particular it follows from Scott's theorem that surface
groups are subgroup separable \cite{scott}. 

In more detail, given
a finite generating set $X$ for $G\subset\S_1$
 and given $h\in\S_1$, one
knows by Scott's theorem that if $h\notin G$ then there
is a finite quotient $\pi:\S_1\to Q$ of $\S_1$ such that 
$\pi(h)\notin\pi(G)$; that is, {\em{$\pi$ separates $h$ from $G$}}.
 To determine if $h\in G$ it is
enough to run two simultaneous processes: on the
one hand one enumerates the finite quotients of $\S_1$
and checks to see if each separates $h$ from $G$; on the
other hand one tries to show that $h\in G$ by simply forming
products $h^{-1}w$ where $w$ is a word in the generators $X$
of $G$, testing to see if each is (freely) equal to a product
of conjugates of the defining relations of $\S_1$.

Now, proceeding by induction on $n$, we assume that
there is a solution to the membership problem for each
finitely presented subgroup of a direct product of $n-1$
or fewer free and surface groups.
Let $D=\S_1\times\dots\times\S_n$ and suppose that
$G$ is a finitely presented subdirect product  of $D$. 
Define  $L_i = G\cap \S_i$.

There is no loss of generality in assuming that elements
$h\in D$ are given as words in the generators of the factors,
and thus we write $h=(h_1,\dots,h_n)$. We assume that the
generators of $G$ are given likewise. 

We first deal with the case where some $L_i$ is trivial, say $L_1$.
The projection of
$G$ to $\S_2\times\dots\times\S_n$ is then isomorphic to $G$, so
in particular it is finitely presented and our induction provides
an algorithm that determines if $(h_2,\dots,h_n)$ lies in
this projection. If it does not, then $h\notin G$. If it does,
then enumerating equalities $h^{-1}w=1$ as above we
eventually find a word $w$ in the generators of $G$ so that
$h^{-1}w$ projects to $1\in \S_2\times\dots\times\S_n$.
Since $L_1=G\cap\S_1=\{1\}$, we deduce that 
in this case $h\in G$ if and only if  $h^{-1}w=1$, and the
validity of this equality can be checked because the word
problem is solvable in $D$.

It remains to consider the case where $G$ intersects
each factor non-trivially. Again we are given $h=(h_1,\dots,h_n)$.
The projection $G_i$ of $G$ to $\S_i$ is finitely generated
and the $\S_i$ are subgroup separable, so we can determine
algorithmically if $h_i\in G_i$. If $h_i\notin G_i$ for some $i$
then $h\notin G$ and we stop. Otherwise, we replace $D$
by the direct product of the $G_i$. Lemma \ref{l:makeP}
allows us to compute a finite presentation for $G_i$
and hence $D$.

We are now reduced to the case where $G$ is a subdirect
product of $D$ and all of the intersections $L_i$
are non-trivial. Again,
Theorem \ref{virtnilp}  tells us that $Q = D/L$ is virtually nilpotent, where 
$L=L_1\times\cdots\times L_n$. Let $\phi:D\to Q$ be the 
quotient map.

Virtually nilpotent groups are subgroup
separable, so if $\phi(h)\notin \phi(G)$ then there
is a finite quotient of $Q$ (and hence $D$) that
separates $h$ from $G$. But $\phi(h)\notin \phi(G)$
if $h\notin G$ because  $L=\ker \phi$ is contained in $G$.
Thus, as in the second paragraph of the proof, an
enumeration of the
finite quotients of $D$ provides an effective procedure
for proving that $h\notin G$ if this is the case.
(Note that we need a finite presentation of $D$ in order
to make this enumeration procedure effective; hence our
earlier invocation of Lemma \ref{l:makeP}.) 

We now have
a procedure that will terminate in a proof if $h\notin G$.  
Once again, we  run this procedure 
 in parallel with a simple-minded enumeration of $h^{-1}w$ that will
 terminate with a proof that $h\in G$ if this is true. 
\end{proof}

\begin{remark} One would like to strengthen the above statement
and claim that there is an algorithm which,
given any  finite set of elements $h, g_1,\ldots,g_k$ of $D$
with  $G =\< g_1,\ldots,g_k\>$  finitely 
presentable, will determine whether or not $h\in G$. But the
above proof fails to establish this precisely
 because the algorithm we
used assumes a knowledge  of which of the intersections
$L_i$ are trivial. 

Our algorithm does provide a uniform algorithm for finitely
presented  products that intersect all of the factors
non-trivially. In fact one can show that such subgroups are all closed
in the profinite topology. With more effort one can show
that the same is true for arbitrary finitely presented subgroups
in a product of 3 surface groups, but
the general situation remains unclear.
\end{remark}

\section{Fibre products and HNN extensions}\label{s:hnn}

Consider an HNN-extension
of the form $$ A = \<B, t\| t^{-1}h t = \phi(h) \ (h\in H)\>$$
where $\phi$ is an isomorphism between a subgroups $H$ and $\phi(H)$ of $B$.
Recall that $A$ is said to be an {\em ascending HNN-extension} if either
$H= B$ or $\phi(H) =B$. In this case either $t$ or $t^{-1}$ conjugates $B$ into
a subgroup of itself.

The metabelian Baumslag-Solitar groups ${\rm{BS}}(1,p) = \<b, t \| t^{-1}bt = b^p\>$
are examples of ascending HNN-extensions.  In \cite{bbhm} it is shown,
for instance, that the untwisted fibre product of two copies of ${\rm{BS}}(1,p)$ mapping
onto the infinte cycle $Q=\<t\|\>$ is finitely presented.

The group ${\rm{BS}}(2,3) = \<b, t \| t^{-1}b^2t = b^3\>$ is a
non-ascending HNN-extension of $\<b\|\>$.  We are going to show, in contrast to \cite{bbhm}, that
the untwisted fibre product $G$ of two copies of ${\rm{BS}}(2,3)$ mapping onto $Q$
is finitely generated but not finitely presented.  Interestingly
we can also show $H_2(G,\Z) = 0$, so that homology does not
detect the lack of a sufficient finite set of relators.

The following gives a large collection of untwisted
fibre products which are finitely generated but not finitely presented. 

\begin{thm}  Let $A_1$ and $A_2$ be non-ascending HNN extensions
with finitely presented base groups and finitely generated 
amalgamated subgroups and
stable letters $t_1$ and $t_2$.  Let $q_i:A_i\to Q=\< t \mid\ \>$ be the map
defined by sending $t_i$ to $t$ and the base groups to 1.  Then the
fibre product $G$ of $q_1$ and $q_2$ is finitely generated but not 
finitely presented.
\end{thm}

\begin{proof}  To simplify notation, even though they are not assumed
to be isomorphic,  we suppress the subscripts
on $A_1$ and $A_2$ for the first part of the proof, adding subscripts
when we consider the fibre product.  By hypothesis we have
a non-ascending HNN extension  
$ A = \<B, t\| t^{-1}h t = \phi(h) \ (h\in H)\>$
where $B$ is finitely presented and $H$ is finitely generated.
Of course $A$ is  finitely presented as well. 

It is easy to see that
the kernel $L$ of $\phi$ has the structure of a two-way infinite, 
proper amalgamated free product.  Setting $B^{i} = t^{-i}Bt^i$,
we observe that $B^0$ and $B^1$ generate their amalagamated free
product 
$$\<B^0, B^1\>= \<B, t^{-1}Bt\> = B 
\underset{\phi(H) = t^{-1}Ht}{\star} (t^{-1}Bt).$$    
Conjugating by $t$ translates this amalgamation decomposition and
$L$ has an amalgamation decomposition indexed by the integers as
$$L =  \cdots \star B^{-1} \underset{t\phi(H)t^{-1} = H}{\star} B^0   
\underset{\phi(H) = t^{-1}Ht}{\star} B^1 \star \cdots.$$
Notice that each amalgamated subgroup is properly contained
in each factor because of our assumption that $A$ is not
an ascending HNN-extension. It follows that $L$ is not finitely
generated.

At this point we reinstate the subscripts on the various objects
associated with $A_1$ and $A_2$.

Now the pullback $G\subset A_1\times A_2$
is generated by  $B_1$ and $B_2$  
together with   $\vc{t}=(t_1,t_2)$.
Also $G\cap A_i = L_i$, the kernel of the map onto
the infinite cycle, has a decomposition as above.  Furthermore,
$G$ is the split extensison of $L_1\times L_2$ by
the infinite cyclic group generated by $\vc{t}=(t_1,t_2)$.

A presentation
for $G$ can be obtained by taking as generators
$\vc{t}$ together with generators for $B_1$ and $B_2$, and
taking as relations:
\begin{enumerate}
\item the finite sets of relations of both $B_1$ and $B_2$;
\item the relations $\vc{t}^{-1}h_1 \vc{t} = \phi_1(h_1)$
and  $\vc{t}^{-1}h_2 \vc{t} = \phi_2(h_2)$  where $h_1$ ranges
over a finite set of generators for $H_1$ and $h_2$ ranges
over a finite set of generators for $H_2$;
\item the relations $u_1v_2 = v_2 u_1$ where $u_1$ ranges over
a set of generators for $L_1$ and $v_2$ ranges over a set
of generators for $L_2$.
\end{enumerate}
The relations in (1) and (2) are finite in number, but the 
relations in (3) are necessarily infinite in number since
neither $L_1$ nor $L_2$ is finitely generated.

Now, in order to obtain a contradiction,
 assume  that $G$ is finitely
presented.  Then in the above presentation some finite
subset $S$ of the relations in (3) together with (1) and (2)
suffice to present $G$.  We can
assume there is a finite portion of the decomposition
of each  $L_i$ of the form
$$B_i^{(m,n)} = B_i^m \underset{(\phi(H))^m = H^{m+1}}{\star}\cdots 
\underset{(\phi(H))^{n-1}= H^n}{\star} B_i^n$$ 
so that all the generators that appear in $S$ lie in
either $B^{(m,n)}_1$ or $B^{(m,n)}_2$. We choose the interval
$(m,n)$ large enough to contain both 0 and 1. Of course each $B_i^{(m,n)}$
is finitely presented.  We are free to add
finitely many relations to those already present, 
so we now enlarge $S$ to contain the finitely many relations
required to say that $B^{(m,n)}_1$ and $B^{(m,n)}_2$
commute.

Observe that the group presented in this way can also be described
as the HNN-extension of  $B^{(m,n)}_1\times B^{(m,n)}_2$
by the stable letter $\vc{t}$ which acts in the same way on
each direct factor and is non-ascending on each.  But then
the normal closure of $B^{(m,n)}_1\times B^{(m,n)}_2$ has a
two-way infinite amalgamated free product decomposition of the
form
$$\cdots\star (B^{(m,n)}_1\times B^{(m,n)}_2)
\underset{M^{(m,n)}}{\star} 
 (B^{(m+1,n+1)}_1\times B^{(m+1,n+1)}_2) \star\cdots$$
Now by our non-ascending assumption there are elements
$x_1 \in B^{(m,n)}_1$ and $y_2\in B^{(m+1,n+1)}_2$
which do not lie in the amalgamation $M^{(m,n)}$.
Hence $x_1y_2x_1^{-1}y_2^{-1} \neq 1$ or equivalently
$x_1y_2\neq y_2x_1$ in this group.  But these elements
clearly must commute in $G$ which is a contradiction.
Hence $G$ could not have been finitely presented.
\end{proof}

\smallskip 

The following gives another example of the type given by Baumslag 
in \cite{gbmult}.

\begin{thm} \label{nfph20} Let $ A = {\rm{BS}}(2,3) = \<b, t \| t^{-1}b^2t = b^3\>$
and let $q: A\to Q = \<t\|\>$ be the map
defined by sending $t$ to $t$ and $b$ to 1.
Then the untwisted fibre product $G$ of two copies of 
the map $q$ is finitely generated but not finitely presented,  and
has $H_2(G,\Z)=0$.
\end{thm}

\proof  An easy calculation shows that $H_1(A,\Z) = \Z$ and  
by the Mayer-Vietoris sequence for HNN-extensions (or the homology
theory of one-relator groups) we have $H_2(A,\Z) = 0$.  
Define $L= [A,A]$ the derived group of $A$.   Then $A$ is also a split
extension of $L$ by $Q$.   The $E^2$ term of the spectral sequence
for this extension has zero maps as differentials so $E^2=E^{\infty}$.  
Since $H_1(A,\Z) = \Z$,
it follows from the spectral sequence that $H_0(Q, H_1(L,\Z)) = 0$.
Also since $H_2(A,\Z) = 0$ it follows that $H_1(Q, H_1(L,\Z)) = 0$
and $H_0(Q, H_2(L,\Z)) = 0$.

Now $L$ has a decomposition as a two-way infinite amalgamated free
product.   If we put $b_i= t^{-i}bt$ for $i\in \Z$ then $b^2_{i+1} = b^3_i$
and the decomposition of $L$ is
$$\cdots \star \<b_{-1}, b_0\mid b^3_{-1}=b^2_0\> \underset{\langle
b_0\rangle}{\star}
 \<b_0, b_1\mid b^3_0=b^2_1\> \underset{\langle
b_1\rangle}{\star}
   \<b_1, b_2\mid b^3_1=b^2_2\>\cdots$$
Computing $H_1(L,\Z)$ by abelianizing this group, one can show that
$H_1(L,\Z)$ is locally cyclic.  Using the (abelian) calculation
$(b_1b_{-1}b_0^2)^6 = b_1^6 b_{-1}^6 b_0^{-12}
=b_0^9b_0^4b_0^{-12} = b_0$
it follows that $H_1(L,\Z) \isom  \Z[\frac{1}{6}]$ and the action induced by
conjugation by $t$ (respectively $t^{-1}$) on  $H_1(L,\Z)$ 
is multiplication by $\frac{3}{2}$ (respectively $\frac{2}{3}$).  

Let $K_{m,n}$ be the subgroup of $L$ generated by $\{b_m,\ldots , b_n\}$.
Again an easy Mayer-Vietoris sequence calculation shows that 
$H_1(K_{0,1},\Z)  = \Z$ and $H_2(K_{0,1},\Z) =0$ (it is the group of the
trefoil knot).  Using the amalgamated free product decomposition
$K_{m,n+1} = K_{m,n} \underset{b_n=b_n}{\star}K_{n,n+1}$,  the 
Mayer-Vietoris sequence shows inductively that 
$H_1(K_{m,n},\Z)  = \Z$ and $H_2(K_{m,n},\Z) =0$.  Since homology
commutes with  direct limits of groups, it follows that $H_2(L,\Z)=0$.

Let $A_1$ and $A_2$ be two copies of $A$ and view $G$ as
a subgroup of $A_1\times A_2$.  Clearly $L_i = G\cap A_i$ and so
$G$ is the split extension of $L_1\times L_2$ by the infinite cycle
generated by $\vc{t} = (t_1,t_2)$ which we can identify with $Q$.
The previous theorem shows that $G$ is finitely generated but
not finitely presented.

Now $H_1(L_1\times L_2,\Z) = H_1(L_1,\Z)\oplus H_1(L_2,\Z) 
= \Z[\frac{1}{6}]\oplus \Z[\frac{1}{6}]$.  Since $H_2(L,\Z)=0$ the Kunneth
theorem shows  $H_2(L_1\times L_2,\Z) \isom H_1(L_1,\Z)\otimes H_1(L_2,\Z)$ 
where conjugation by $\vc{t}$ induces the diagonal action on the
right hand tensor product.  Since
$H_1(L,\Z) = \Z[\frac{1}{6}]$,  the calculation 
$\frac{1}{6^n}\otimes 1= \frac{1}{6^n}\otimes \frac{6^n}{6^n} 
= \frac{6^n}{6^n}\otimes \frac{1}{6^n}
=1 \otimes \frac{1}{6^n}$ implies
$H_1(L_1,\Z)\otimes H_1(L_2,\Z)\isom H_1(L,\Z)$.  Hence we
have $H_2(L_1\times L_2,\Z) \isom H_1(L,\Z)$.

We now examine the $E^2$ terms $H_p(Q, H_q(L_1\times L_2,\Z))$
of the spectral sequence for $G$ with $p+q=2$.  Since $Q$ is
the infinite cyclic group  $H_2(Q, H_0(L_1\times L_2,\Z))=0$.  
Using the early results about the spectral sequence for $A$ we have
$$H_1(Q, H_1(L_1\times L_2,\Z) \isom H_1(Q, H_1(L_1,\Z)\oplus H_1(L_2,\Z)$$
$$\hfill \isom H_1(Q, H_1(L_1,\Z)) \oplus H_1(Q, H_1(L_1,\Z))  = 0\oplus 0 = 0$$
and
$$H_0(Q,H_2(L_1\times L_2,\Z)) \isom H_0(Q, H_1(L,\Z)) = 0.$$
Since all three of these terms with  $p+q=2$ are 0, it follows
that $H_2(G,\Z) = 0$.  This completes the proof. 
\qed

\smallskip

\noindent{\bf{Acknowledgements:}} The following organizations provided financial
support at various stages of this project:  the Australian Research Council,
the EPSRC of the United Kingdom,
 the Swiss National
Science Foundation and l'Alliance Scientific (grant \# PN 05.004).
Bridson is also supported in part by a 
Royal Society Nuffield Research Merit Award. We are grateful
to all of these organizations.

We also thank l'Universit\'e
de Provence Aix-Marseilles, the EPFL (Lausanne) and l'Universit\'e
de Gen\`eve for their hospitality during the writing of this article,
and we particularly thank Hamish Short for his extended hospitality.

\end{document}